\documentclass[12pt,twoside]{amsart}\usepackage{amssymb}
\usepackage{enumerate}

\newtheorem{thmspec}{\relax}

\newtheorem{theorem}{Theorem}[section]
\newtheorem{thm}[theorem]{Theorem}

\newtheorem{cor}[theorem]{Corollary}
\newtheorem{prop}[theorem]{Proposition}

\newtheorem{rems}[theorem]{Remarks}
\theoremstyle{definition}

\theoremstyle{remark}

\numberwithin{equation}{section}
\def \cal{\mathcal}
\def \Bbb{\mathbb}

\def\onto{{\kern3pt\to\kern-8pt\to\kern3pt}}

\def\<{\langle}
\def\>{\rangle}
\def\|{{\ |\ }}

\def\onto{\twoheadrightarrow}

\def\-{\underline}

\def\N{\Bbb N}

\def\R{\Bbb R}
\def\C{\Bbb C}

\def\L{\cal L}

\def\dim{\operatorname{dim}}

\def\d{\operatorname{d}}
\def\Ker{\operatorname{Ker}}

\def\ddc{\operatorname{dd^c}}
\def\Vol{\operatorname{Vol}}

\def\P{\mathcal{P}}

\def\H{\mathcal{H}}
\def\K{\mathcal{K}}
\def\E{\mathcal{E}}

\def\<{\langle}
\def\>{\rangle}

\catcode`\@=11
\def\serieslogo@{\relax}
\def\@setcopyright{\relax}
\catcode`\@=12

\title[The mixed Hodge--Riemann bilinear relations]
{The mixed Hodge--Riemann bilinear relations\\
for  compact  K\"{a}hler manifolds}

\begin{document}

\author{Tien-Cuong  Dinh}
\address{Tien-Cuong  Dinh \\
  Math\'ematique-B\^at 425\\
 UMR 8628, Universit\'e Paris-Sud
  F--91405 Orsay, France}
\email{TienCuong.Dinh@math.u-psud.fr}

\author{Vi\^et-Anh  Nguy\^en}
\address{Vi\^et-Anh  Nguy\^en\\
 Max-Planck Institut f\"{u}r    Mathematik\\
Vivatsgasse 7, D--53111\\
Bonn, Germany}
\email{vietanh@mpim-bonn.mpg.de}



\maketitle

\begin{center}
{\small \it Dedicated to Professor Henri Skoda on the occasion of his 60th
birthday} 
\end{center}

\begin{abstract}
We prove  the   Hodge--Riemann bilinear relations, the hard Lefschetz 
theorem and the Lefschetz decomposition 
for  compact  K\"{a}hler manifolds in the mixed situation. 
 \end{abstract}

\section{Introduction  and statement of the main results}

Around the year 1979,   Khovanskii (see \cite{ko1,ko2,ko}) and Teissier (see \cite{te1,te,te2})  discovered independently  a beautiful   intimate relationship between
the theory of mixed volumes and algebraic geometry.  In order to describe  this connection  we recall some facts from the theory of mixed volumes.
Let $K_1,\ldots, K_r$  be $r$  $n$-dimensional convex bodies  in $\R^n$ and $I=(i_1,\ldots,i_r)\in \N^r$ with $\vert I\vert:=\sum\limits_{s=1}^r i_s=n.$
Then  the (Minkowski) mixed volume $K^I=\left [K_1^{i_1} \ldots K_r^{i_r} \right ] $  is determined by the folowing  identity
\begin{equation*}
\Vol\Big(\sum\limits_{s=1}^{r}\lambda_s K_s\Big)=\sum\limits_{I=(i_1,\ldots,i_r):\ \vert I\vert=n}\frac{n!}{i_1!\cdots i_r!} K^I \lambda_1^{i_1}\cdots\lambda_r^{i_r}, \qquad
\lambda_1,\ldots,\lambda_r\geq 0.
\end{equation*}
The  Aleksandrov--Fenchel  inequalities state that
\begin{equation*}
\Big(\left [K_1 K_2 \ldots K_r \right ] \Big)^2 \geq \left [K_1 K_1 K_3 \ldots K_r \right ]\cdot \left [K_2 K_2 K_3 \ldots K_r \right ].
\end{equation*}
    
    Now let $X$ be a   complex  algebraic manifold of  dimension $n$
  and  $D_1,\ldots,D_r$  very  ample divisors on $X.$ Let $D^I=\left [D_1^{i_1} \ldots D_r^{i_r} \right ] $  denote the index of intersection of 
 $ D_1^{i_1}\cap\cdots\cap  D_r^{i_r}  ,$ where $D_s^{r_s}$ stands for $D_s\cap\cdots\cap D_s$  ($r_s$ times).
  Khovanskii   and Teissier found out a profound  analog between   Aleksandrov--Fenchel  inequalities  and the Hodge--Riemann bilinear relations   in algebraic geometry
 \begin{equation*}
\Big(\left [D_1 D_2 \ldots D_r \right ]\Big)^2 \geq \left [D_1 D_1 D_3 \ldots D_r \right ]\cdot \left [D_2 D_2 D_3 \ldots D_r \right ].
\end{equation*} 
  Their proofs use the usual Hodge--Riemann bilinear relations (see Theorem \ref{classicHR} below) applied to K\"{a}hler forms corresponding to the divisors  
and an induction argument. 
  Khovanskii   and Teissier also noted that many other  interesting inequalities from convex geometry (for example the Brunn--Minkowski inequality,
  Bonnesen-type inequalities etc.) either could
  be deduced from  the Hodge--Riemann bilinear relations, or  find  their analogs for algebraic varieties that 
  generalize the   Hodge--Riemann bilinear relations. Based on  this point of view
  P. McMullen (see \cite{mm}) developed   a deep and important generalization of Aleksandrov--Fenchel  inequalities for simple convex polytopes. 
  On the other hand, Khovanskii   and Teissier's discovery   also suggests a generalization of
  the mixed Hodge--Riemann bilinear relations in the context of compact K\"ahler manifolds.
That is the main motivation of our work.

Let $X$ be a  compact K\"ahler manifold of dimension $n$.
 Let  $0\leq p,q\leq n$ and $0\leq r\leq 2n$ be integers.  
One denotes by $\E^{p,q}(X)$  (resp.  $L^2_{p,q}(X)$)   the  space of complex-valued 
   differential forms of bidegree $(p,q)$ on $X$ with smooth 
coefficients (resp. with $L^2$-coefficients).
   For  $\alpha\in L^2_{p,q}(X),$ $\Vert\alpha\Vert_{L^2}$ denotes  
its $L^2$-norm, i.e. the sum of $L^2$-norms of its coefficients
   on charts.  In the sequel 
   $\H^{p,q}(X)$ denote  the space of smooth $\d$-closed  $(p,q)$-forms modulo
    smooth $\d$-exact  $(p,q)$-forms.
     Moreover, for any smooth $\d$-closed form 
$\alpha\in \E^{p,q}(X),$ 
   $[\alpha]$ denotes the class  of $\alpha$ in  
$\H^{p,q}(X) .$  We can identify $\H^{p,q}(X)$ to the subspace
of $\H^{p+q}(X)$ spanned by classes of  smooth $\d$-closed  $(p,q)$-forms.
The classical Hodge decomposition theorem asserts that
\begin{equation*}
\H^{r}(X)=\bigoplus\limits_{p+q=r} \H^{p,q}(X)
\qquad\text{and}\qquad  \H^{p,q}(X)=\overline{\H^{q,p}(X)}.
\end{equation*}
We refer the reader to \cite{bdi,gh,vo,we}   for the basics of 
the Hodge theory and to
\cite{bs,gr,ko,te, ti} for some of its advanced aspects.

Fix  non-negative integers  $p,\ q$  such that $p+q\leq n.$  Let $\omega_1,\ldots,\omega_{n-p-q+1}$ be  K\"{a}hler forms.
Put $\Omega:=\omega_1\wedge\cdots\wedge\omega_{n-p-q}$.
Consider the {\it mixed primitive subspace}
   \begin{equation}\label{eq1.1}
   P^{p,q}(X):= \left \lbrace  [\alpha]\in \H^{p,q}(X) :\ [\alpha] \wedge[\Omega]\wedge [\omega_{n-p-q+1}]=0 \right\rbrace.
   \end{equation}
 Let us define the  {\it mixed Hodge--Riemann   bilinear  form} on $\H^{p,q}(X)\otimes \H^{p,q}(X)$ as follows
 \begin{equation}\label{eq1.2}
Q([\alpha],[\beta]):=i^{p-q}(-1)^{\frac{(n-p-q)(n-p-q-1)}{2}}\int_{X} 
\alpha\wedge \overline{\beta}\wedge  \Omega .
\end{equation} 
Observe that  $Q(\cdot,\cdot)$ is a sesquilinear Hermitian symmetric form.

\smallskip\smallskip

The classical  Hodge--Riemann bilinear relations 
state that
\begin{thm}\label{classicHR}
If  $\omega_1 =\cdots=\omega_{n-p-q+1},$  then  $Q(\cdot,\cdot)$ is positive definite on the primitive space
$P^{p,q}(X).$
\end{thm}

The open question  can be formulated as follows:

{\it  Does  Theorem  \ref{classicHR} still  hold if $\omega_1,\ldots,\omega_{n-p-q+1}$ are arbitrary K\"{a}hler forms?  }
 
An attempt towards this generalization is 
made by Gromov. Namely,    the following theorem is stated in  \cite{gr}.

\begin{thm}  (Gromov's Theorem)
  If $p=q,$  
  then  $Q(\cdot,\cdot)$ is  positive semi-definite on  
$P^{p,q}(X)$, that is, $Q([\alpha],[\alpha])\geq 0$  for $\alpha\in P^{p,q}(X)$. 
\end{thm}
However, Gromov only gave therein a 
complete proof for the special case where $p=q=1.$ 
On continuation of Gromov's work and using Aleksandrov's approach,
Timorin has  proved a general mixed Hodge--Riemann bilinear relations,  but only in the linear situation \cite{ti} (see also \cite{ko, te}).
His result may be rephrased as follows (see also Proposition \ref{timorin1} below).
\begin{thm}  (Timorin's Theorem)\label{timorin}
  If $ X$ is a complex torus of dimension $n,$ 
  then  $Q(\cdot,\cdot)$ is positive definite on  
$P^{p,q}(X).$  
\end{thm}

The purpose of this article is to prove the above theorems in the general context.
Now we state the  main results.
\renewcommand{\thethmspec}{Theorem A}
  \begin{thmspec} Let $X$ be a  compact K\"{a}hler manifold   of dimension $n$ and    $p,q$   integers such that $0\leq p,q\leq p+q\leq n.$ Then,
for arbitrary K\"{a}hler forms  $\omega_1,\ldots,\omega_{n-p-q+1},$ the mixed Hodge--Riemann   bilinear  form $Q(\cdot,\cdot)$ is positive definite on the mixed primitive subspace
$P^{p,q}(X).$
\end{thmspec}

Note that when $\omega_j$  are cohomologous to very ample divisors
of $X,$ by Bertini theorem, one can replace $[\omega_j]$ by divisors $D_j$ 
which
intersects transversally. Then one deduces   from the
classical Hodge-Riemann theorem on the submanifold $D:=D_1\cap \cdots \cap 
D_{n-p-q}$ that $Q([\alpha],[\alpha])\geq 0$  for 
all $[\alpha]  $ satisfying  $[\alpha]\wedge [\omega_{n-p-q+1}]=0$ on $\H^{p+1,q+1}(D)$
(see also \cite{ko1,ko2,te} and \cite{vo}). This is the original reason
to believe that the mixed  Hodge--Riemann bilinear relations hold in the general situation.

The  following results  
generalize the hard Lefschetz theorem and the Lefschetz decomposition theorem. 
\renewcommand{\thethmspec}{Theorem B}
  \begin{thmspec}  Let $X$ be a  compact K\"{a}hler manifold   of dimension $n$ and    $p,q$   integers such that $0\leq p,q\leq p+q\leq n.$ 
  Then,
for arbitrary K\"{a}hler forms  $\omega_1,\ldots,\omega_{n-p-q},$ the  linear map $\tau:\ \H^{p,q}(X)\longrightarrow \H^{n-q,n-p}(X)$ given by
\begin{equation*}
\tau([\alpha]):=[\Omega]\wedge[\alpha],\qquad  [\alpha]\in \H^{p,q}(X),
\end{equation*}
where  $[\Omega]:=[\omega_1]\wedge\cdots\wedge[\omega_{n-p-q}],$
is an isomorphism.
\end{thmspec}

\renewcommand{\thethmspec}{Theorem C}
  \begin{thmspec}  Let $X$ be a  compact K\"{a}hler manifold   
of dimension $n$   and    $p,q$   integers such that 
$0\leq p,q\leq p+q\leq n.$  Then,
for arbitrary K\"{a}hler forms  $\omega_1,\ldots,\omega_{n-p-q+1},$ 
the  following canonical decomposition holds
\begin{equation*}
 \H^{p,q}(X)=P^{p,q}(X)\oplus [\omega_{n-p-q+1}]\wedge \H^{p-1,q-1}(X),
\end{equation*}
with the convention that  $\H^{p-1,q-1}(X):=0$ if either $p=0$ or $q=0.$  
 \end{thmspec}

  We close  the introduction with a brief outline of the paper to follow.
 
 \smallskip
 
  Our strategy is to reduce the general case to the linear case.  
   In order to achieve this reduction  we apply the $L^2$-technique 
to solve a $\ddc$-equation. Recall here that ${\rm d}=\partial +\overline\partial$,
${\rm d}^{\rm c}=\frac{i}{2\pi}(\overline\partial-\partial)$ and 
$\ddc=\frac{i}{\pi}\partial\overline\partial$.
   Section 2 is then  devoted to developing the necessary technique.
     We begin  this section by  collecting some results 
of Timorin and by establishing
  some estimates.  This will enable us to construct a solution 
of  the above equation.  We will, in the remaining part of Section 2,
  regularize this solution.  
Based on the results of Section 2, the proofs of the main 
theorems are presented in Section 3.
  
The mixed Hodge-Riemann theorem is not true in general if we replace 
$[\Omega]$ by the class of a smooth strictly positive form 
as a simple example in \cite{bs}
shows. However,
by continuity, it holds for every class $[\Omega]$ close enough to a product 
of K\"ahler classes. In Section 4 we describe the domain of validity 
of this theorem in the case where $p=q=1$.

  \smallskip

\noindent{\bf Acknowledgment.}    We would like to thank the referee for many
interesting   suggestions and remarks. We are also grateful to Professor Nessim Sibony for very stimulating discussions. 
The second author  wishes to express his gratitude to the Max-Planck Institut f\"{u}r Mathematik
in Bonn (Germany)  for its hospitality and its support.

\section{Preparatory results}
In the first two propositions we  place ourselves in the linear context. 
For $0\leq p,q\leq n,$  let $\Lambda^{p,q}(\C^n)$  
denote the space of $(p,q)$-forms with complex-constant coefficients. $\Lambda^{p,q}(\C^n)$ is equipped
with the Euclidean norm $\Vert\cdot\Vert.$  
We first recall  Timorin's result \cite{ti}.
\begin{prop}\label{timorin1} Let $p,q$ be  integers such that 
$0\leq p,q\leq p+q\leq n$  and 
  $\omega_1,\ldots,\omega_{n-p-q+1}$   strictly positive 
forms of $\Lambda^{1,1}(\C^n).$ 
Define the   sesquilinear Hermitian symmetric form
   \begin{equation*} 
Q(\alpha,\beta):=i^{p-q}(-1)^{\frac{(n-p-q)(n-p-q-1)}{2}} \ast\Big(\alpha\wedge \overline{\beta}\wedge  \Omega \Big),\qquad\alpha,\beta\in \Lambda^{p,q}(\C^n),
\end{equation*} 
where  $\ast$ is the Hodge star operator, and 
$\Omega:=\omega_1\wedge\cdots\wedge\omega_{n-p-q}$. 
Define the mixed primitive subspace
  \begin{equation*} 
   P^{p,q}(\C^n):= \left \lbrace  \alpha\in  \Lambda^{p,q}(\C^n):\  \alpha\wedge\Omega\wedge \omega_{n-p-q+1} =0 \right\rbrace.
   \end{equation*} Then
   \begin{enumerate}
\item[{\rm (a)}] The operator of multiplication by $\Omega$  
induces an isomorphism between $\Lambda^{p,q}(\C^n)$ and 
 $\Lambda^{n-q,n-p}(\C^n).$ 
 \item[{\rm (b)}] 
   $Q(\cdot,\cdot)$ is positive definite on  
$P^{p,q}(\C^n).$
 \item[{\rm (c)}]   The space   $\Lambda^{p,q}(\C^n)$ splits into the $Q$-orthogonal direct sum
 \begin{equation*}
 \Lambda^{p,q}(\C^n)=P^{p,q}(\C^n)\oplus \omega_{n-p-q+1}\wedge 
\Lambda^{p-1,q-1}(\C^n),
 \end{equation*}
with the convention that $\Lambda^{p-1,q-1}(\C^n):=0$ if either $p=0$ or $q=0$. 
 \end{enumerate}
\end{prop}
\begin{proof}  See Proposition 1, the Main Theorem  and Corollary 2 in \cite{ti}.
\end{proof}
The following  estimate will be  crucial later on. 
\begin{prop}\label{timorin2} There are finite positive constants $C_1$ and $C_2$ such that
\begin{equation*}
  C_1\cdot\Vert \alpha\wedge\Omega\wedge \omega_{n-p-q+1}\Vert^2 +C_2\cdot
  \Re Q(\alpha,\alpha)\geq  \Vert\alpha\Vert^2
 \end{equation*}
for all forms $\alpha\in  \Lambda^{p,q}(\C^n).$
\end{prop}
\begin{proof} 
  By  Proposition \ref{timorin1}(a) we may find  
a positive finite constant  $C$ so that
 \begin{equation}\label{timorin2.0}
   \frac{\Vert \gamma\Vert}{C} \leq   \Vert \gamma\wedge\Omega\wedge \omega_{n-p-q+1}^2\Vert
\leq C\cdot \Vert \gamma\Vert , \qquad \gamma\in \Lambda^{p-1,q-1}(\C^n).
\end{equation}
 By   Proposition \ref{timorin1}(c) we may write
\begin{equation*}
 \alpha=\beta +  \omega_{n-p-q+1}\wedge\gamma, \qquad \beta\in P^{p,q}(\C^n),\  \gamma\in\Lambda^{p-1,q-1}(\C^n).
 \end{equation*}
 Then   we have
 \begin{equation}\label{timorin2.1}
Q( \alpha,\alpha)=Q(\beta,\beta)+Q(   \omega_{n-p-q+1}\wedge\gamma, \omega_{n-p-q+1}\wedge\gamma).
 \end{equation}
 On the other hand, since  $\beta\in P^{p,q}(\C^n),$ one   gets that
 \begin{equation}\label{timorin2.2}
   \Vert \alpha\wedge\Omega\wedge \omega_{n-p-q+1}\Vert  = \Vert \gamma\wedge\Omega\wedge \omega_{n-p-q+1}^2\Vert
\geq \frac{\Vert \gamma\Vert}{C}, \end{equation}
where the   estimate follows from  the left-side estimate 
in (\ref{timorin2.0}). Therefore, we obtain, for $C'>0$ large enough,
\begin{equation}\label{timorin2.3}
\Vert \alpha\Vert^2\leq C'(\Vert\beta\Vert^2+\Vert\gamma\Vert^2)\leq 
C'\Vert\beta\Vert^2+
   C'C^2\Vert \alpha\wedge\Omega\wedge \omega_{n-p-q+1}\Vert^2.    
   \end{equation}
   On the other hand,  by   Proposition  \ref{timorin1}(b) we may find a 
positive finite constant  $C''$ so that
   \begin{eqnarray*}
   \Vert \beta\Vert^2\leq C''\cdot Q(\beta,\beta)&=& 
C''\cdot\Big (Q(\alpha,\alpha)- Q(   \omega_{n-p-q+1}\wedge\gamma, 
\omega_{n-p-q+1}\wedge\gamma)   \Big) \\
   &=& C''\cdot\Big (\Re Q(\alpha,\alpha)-
\Re Q(   \omega_{n-p-q+1}\wedge\gamma, \omega_{n-p-q+1}\wedge\gamma)  \Big)\\
  & \leq& C''\cdot\Re Q(\alpha,\alpha)+  C''C^2\cdot\Vert \gamma\Vert^2\\
   &\leq& C''\cdot \Re Q(\alpha,\alpha)+ 
C''C^4\cdot \Vert \alpha\wedge\Omega\wedge \omega_{n-p-q+1}\Vert^2,
   \end{eqnarray*}
where the first identity follows from (\ref{timorin2.1}), the second estimate from  the right-side estimate in (\ref{timorin2.0}),  and  the last one from (\ref{timorin2.2}).
This, combined with (\ref{timorin2.3}), implies the desired estimate for 
$C_1:=C'C''$ and $C_2:=C'C''C^4+C'C^2.$  
\end{proof}

\begin{prop}\label{ddc_equation}
 We keep the hypothesis and the notation in the statement of Theorem A.
Assume that $p\geq 1$ and $q\geq 1$.
Then, for every $\d$-closed form $f\in \E^{p,q}(X)$ such that  $ [f]\in P^{p,q}(X),$
there is a form $u\in L^2_{p-1,q-1}(X)$ such that  
\begin{equation*}
\ddc u\wedge\Omega\wedge\omega_{n-p-q+1}=f\wedge\Omega\wedge\omega_{n-p-q+1}.
\end{equation*}
\end{prop}
\begin{proof}
Consider the subspace $H$ of $L^2_{n-p+1,n-q+1}(X)$  defined by
\begin{equation*}
 H:=\left\lbrace  \ddc \alpha\wedge\Omega\wedge\omega_{n-p-q+1}:\  \alpha\in \E^{q-1,p-1}(X)\right\rbrace.
\end{equation*}
We  construct a linear form $h$ on $H$ as follows
\begin{equation}\label{ddc_2}
  h\left( \ddc \alpha\wedge\Omega\wedge\omega_{n-p-q+1}\right):= (-1)^{p+q}\int_{X}\alpha\wedge f\wedge  \Omega\wedge\omega_{n-p-q+1}. 
\end{equation}
We now check that  $h$ is a well-defined bounded linear form with respect to the $L^2$-norm restricted to $H.$
 To this end one first shows that there is a positive finite constant $C$ such that
 \begin{equation}\label{ddc_3}
    \left \Vert \ddc \alpha\right\Vert_{L^2}\leq C\cdot \left \Vert \ddc \alpha\wedge\Omega\wedge\omega_{n-p-q+1}\right\Vert_{L^2}  .
\end{equation}
  To prove  (\ref{ddc_3}) we first use a compactness argument to find  finite disjoint open sets $(U_j)_{j=1}^N$ of $X$ so that  $\overline{U_j}$ is    contained in a  local  chart,
   and that $\partial U_j$  is piecewisely smooth, and that $X=\bigcup\limits_{j=1}^N \overline{U_j}.$ One next invokes the estimate in Proposition \ref{timorin2} for every
     point in each  $U_j,$
$j=1,\ldots,N.$  Then one integrates this estimate  over $X.$  
We extend the bilinear form $Q(\cdot,\cdot)$  given by  formula (\ref{eq1.2}) in a canonical way to
  $\E^{p,q}(X)\otimes  \E^{p,q}(X):$ 
\begin{equation*} 
Q( \alpha , \beta ):=i^{p-q}(-1)^{\frac{(n-p-q)(n-p-q-1)}{2}}\int_{X} 
\alpha\wedge \overline{\beta}\wedge  \Omega ,\qquad  \alpha,\beta\in  \E^{p,q}(X).
\end{equation*} 
 Consequently,
 for suitable  positive finite constants  $C$ and $C^{'},$
 \begin{equation}\label{ddc_4} 
   \left \Vert \ddc \alpha\right\Vert_{L^2}^2\leq C\cdot \left \Vert \ddc \alpha\wedge\Omega\wedge\omega_{n-p-q+1}\right\Vert_{L^2}^2  +
   C^{'}\cdot\Re  Q(\ddc \alpha,\ddc\alpha).
\end{equation}
 On the other hand, applying Stokes' Theorem yields that
\begin{equation*} 
      Q(\ddc \alpha,\ddc\alpha)= i^{p-q}(-1)^{\frac{(n-p-q)(n-p-q-1)}{2}}
     \int_{X} \ddc \alpha\wedge\ddc\overline{\alpha}\wedge\Omega=0.
\end{equation*}
 This, combined with  (\ref{ddc_4}), implies   (\ref{ddc_3}).  
 
By hypothesis the smooth form $f\wedge\Omega\wedge\omega_{n-p-q+1}$ is $\d$-exact. 
Consequently, it follows
from \cite[p. 41]{bdi}  that  there is a form $g\in \E^{n-q,n-p}(X)$ such that 
\begin{equation*}
\ddc g =f\wedge\Omega\wedge\omega_{n-p-q+1}.
\end{equation*}
Applying Stokes' Theorem, we obtain that 
 \begin{eqnarray*}
 \left\vert \int_{X}\alpha\wedge f\wedge  \Omega\wedge\omega_{n-p-q+1}\right\vert
& = & 
 \left\vert\int_{X}\alpha\wedge \ddc g\right\vert=\left\vert\int_{X}\ddc \alpha\wedge  g\right\vert\\
 & \leq &   \Vert g\Vert_{L^2}\cdot  \left \Vert \ddc \alpha\right\Vert_{L^2} \\
& \leq &  C \Vert g\Vert_{L^2}\cdot \left \Vert \ddc \alpha\wedge\Omega\wedge\omega_{n-p-q+1}\right\Vert_{L^2},
 \end{eqnarray*}
 where  the latter estimate follows from (\ref{ddc_3}).  In particular, we have 
 \begin{equation*}
 \int_{X}\alpha\wedge f\wedge  \Omega\wedge\omega_{n-p-q+1}= 0 \qquad \text{when}\quad \ddc \alpha\wedge\Omega\wedge\omega_{n-p-q+1} =0.
 \end{equation*}
 In summary, we have just shown that  $h$ given 
by (\ref{ddc_2}) is a well-defined bounded linear form with respect to the $L^2$-norm restricted to $H,$
 and its norm is dominated by $C \Vert g\Vert_{L^2}.$ Applying the Hahn--Banach Theorem, we may extend $h$ to a bounded linear form  on $L^2_{n-p+1,n-q+1}(X).$
 Let $u$ be a form in  $L^2_{p-1,q-1}(X)$  that represents $h.$ Then, 
in virtue of (\ref{ddc_2}), we have that
 \begin{equation*} 
   \int_{X} u\wedge \ddc\alpha \wedge\Omega\wedge\omega_{n-p-q+1}  = (-1)^{p+q}\int_{X}\alpha\wedge f\wedge  \Omega\wedge\omega_{n-p-q+1}
\end{equation*}
for all test forms $\alpha\in \E^{q-1,p-1}(X).$ This is the desired identity of the proposition.
\end{proof}

We need to regularize the solution $u$ given by the previous proposition. This is the purpose of the following result.
\begin{prop}\label{smooth_solution}
 We keep the hypothesis  and the conclusion in the statement of  Proposition \ref{ddc_equation}.
   Then,  
there is a form $v\in \E^{p-1,q-1}(X)$ such that  $\ddc v=\ddc u.$
 \end{prop}
 \begin{proof}
 First we like to equip the vector bundle $\E^{p,q}(X)$ with a special Hermitian metric.
 To this end suppose without loss of generality that  $p\leq q.$
 For any $\alpha\in \E^{p,q}(X),$ we apply  Proposition \ref{timorin1}(c) repeatedly in order to obtain 
 the following unique decomposition
 \begin{equation}\label{smooth_solution_eq1}
 \alpha=\sum\limits_{j=0}^p \alpha_j\wedge \omega_{n-p-q+1}^{p-j},
 \end{equation} 
 where $\alpha_j\in \E^{j,q-p+j}(X)$  such that $\alpha_j\wedge \Omega\wedge \omega_{n-p-q+1}^{2p-2j+1}=0$
 (see also (\ref{eq1.1})).
 Now we can define a new form $\widetilde{\alpha}\in \E^{p,q}(X)$ as follows
 \begin{equation}\label{smooth_solution_eq2}
 \widetilde{\alpha}:=\sum\limits_{j=0}^p (-1)^{p-j}\alpha_j\wedge \omega_{n-p-q+1}^{p-j}.
 \end{equation} 
 Define an inner product $\left <\cdot, \cdot\right>$ on $\E^{p,q}(X)$  by setting
 \begin{equation}\label{smooth_solution_eq3}
 \left <\alpha, \beta\right> := Q(\alpha,\widetilde{\beta} ),\qquad
\alpha, \beta\in \E^{p,q}(X),
 \end{equation} 
 where  $Q(\cdot,\cdot)$ is given by the same integral as in (\ref{eq1.2}). 
Using Proposition \ref{timorin1}(c),
 one may rewrite (\ref{smooth_solution_eq3}) as follows
 \begin{equation*} 
 \left <\alpha, \beta\right>= \sum\limits_{j=0}^p (-1)^{p-j}Q(  \omega_{n-p-q+1}^{p-j}\wedge \alpha_j,\omega_{n-p-q+1}^{p-j}\wedge \beta_j ),
 \end{equation*} 
 where  the  $\beta_j$'s are determined by $\beta$   in virtue of   (\ref{smooth_solution_eq1}).
 Applying    Proposition \ref{timorin1}(b) and using (\ref{eq1.2}) and  (\ref{smooth_solution_eq1})--(\ref{smooth_solution_eq3}), one can check that
$\left <\cdot, \cdot\right>$ defines a Hermitian metric on $\E^{p,q}(X).$ Moreover, if we consider the norm 
$\Vert \alpha\Vert:=\sqrt{\left <\alpha, \alpha\right>},$ then   there is a positive finite constant $C$ such that
 \begin{equation*} 
 \frac{1}{C}\cdot \Big ( \sum\limits_{j=0}^p   \Vert \alpha_j\Vert_{L^2} \Big) \leq  \Vert \alpha\Vert\leq C\cdot  \sum\limits_{j=0}^p   \Vert \alpha_j\Vert_{L^2}.
 \end{equation*}

  Consider the following form of bidegree $(p,q)$
 \begin{equation}\label{smooth_solution_eq4}
  h:=\ddc u-f.
 \end{equation}
 Then in virtue of Proposition \ref{ddc_equation} and of the hypothesis, 
 $h$ belongs to the Sobolev space  $W^{-2}(\E^{p,q}(X))$\footnote{For 
 the Sobolev spaces on compact manifolds, see Chapter IV in \cite{we}}.  In addition, the following identities hold
  \begin{equation}\label{smooth_solution_eq5}
 \overline{\partial} h=0,\quad \partial h=0\quad  \text{and}\quad   h\wedge \Omega\wedge \omega_{n-p-q+1}=0.
 \end{equation}  
 
 For any form $\alpha\in \E^{p,q-1}(X),$ we have that
 \begin{equation*} 
 \left <\overline{\partial}\alpha, h\right>= Q\left(\overline{\partial}\alpha, h\right)=
 i^{p-q}(-1)^{p+q-1+\frac{(n-p-q)(n-p-q-1)}{2}  }\int_{X}\alpha\wedge \overline{\partial}\overline{h}\wedge  \Omega =0,
 \end{equation*} 
 where the first identity follows from (\ref{smooth_solution_eq1})--(\ref{smooth_solution_eq3}) and from the
 third identity in  (\ref{smooth_solution_eq5}), the second one from (\ref{eq1.2}) and from an application of Stokes' Theorem,
 and the last one from the second identity in  (\ref{smooth_solution_eq5}). Let  $\overline{\partial}^{\ast}$ be the adjoint of $\overline{\partial}$
 with respect to the inner product given in (\ref{smooth_solution_eq3}). Then we have shown that $\overline{\partial}^{\ast} h=0.$
 On the other hand,  $\overline{\partial} h=0$ by (\ref{smooth_solution_eq5}) and $h \in W^{-2}(\E^{p,q}(X)).$
 Therefore, $h$ is a harmonic current with respect to the Laplacian operator
 $ \overline{\partial} \overline{\partial}^{\ast}+\overline{\partial}^{\ast}\overline{\partial} $  (see Section 5 in \cite[Chap. IV]{we}).
 Consequently, by elliptic regularity (see Theorem 4.9 in \cite[Chap. IV]{we}) 
 $h$ is smooth. Hence, $\ddc u$ is smooth by (\ref{smooth_solution_eq4}).
 By the classical Hodge theory \cite[p. 41]{bdi}
   there is a $v\in\E^{p-1,q-1}(X)$ such that  $\ddc v=\ddc u.$ 
Hence, the proof is finished. 
 \end{proof}

 \section{Proof of the main results}

 Now we arrive at

\smallskip

\noindent{\bf  Proof of Theorem A.} Let  $f$   be a $\d$-closed form in $ \E^{p,q}(X)$ such that  $ [f]\in P^{p,q}(X).$
We like  to prove that   $Q([f],[f])\geq 0.$ Let $v$ be the smooth $(p-1,q-1)$-form  given by  Proposition \ref{smooth_solution}.  
Then we have
\begin{equation}\label{correction}
(f-\ddc v)\wedge  \Omega\wedge\omega_{n-p-q+1}=0.
\end{equation}
When either $p=0$ or $q=0$ we replace $\ddc v$ by 0.
 In virtue of  the identity (\ref{correction}), 
we are able to apply  Proposition \ref{timorin1}(b) to every point of $X.$
Consequently, after an integration on $X,$  we obtain that  
\begin{equation}\label{TheoremA_eq1} 
i^{p-q}(-1)^{\frac{(n-p-q)(n-p-q-1)}{2}}\int\limits_{X} (f-\ddc v)\wedge (\overline{f}-\ddc \overline{v})\wedge \Omega \geq 0.
\end{equation}
Applying Stokes' Theorem to the left-hand side of the last line yields that
 \begin{equation*} 
 \int\limits_{X} f \wedge \overline{f} \wedge \Omega= \int\limits_{X} (f-\ddc v)\wedge (\overline{f}-\ddc \overline{v})\wedge \Omega  .
\end{equation*}
This, combined with (\ref{TheoremA_eq1}), implies 
  that $Q([f],[f])\geq 0.$ The equality happens if and only if $f=\ddc v,$   in other words,
$[f]=0.$  Hence, the proof of the theorem is complete.  
 \hfill $\square$

\smallskip

\noindent{\bf  Proof of Theorem B.}
Let $\omega_{n-p-q+1}$ be an arbitrary K\"{a}hler form.  Since  $\dim\H^{p,q}(X)=\dim \H^{n-q,n-p}(X),$ it is sufficient to show that $\tau$ is injective.
  To this end let $\alpha$ be a  $\d$-closed form   in $ \E^{p,q}(X)$ such that  
\begin{equation*}
\tau([\alpha])= [\alpha]\wedge [\Omega]=0\ \text{in}\  \H^{n-q,n-p}(X) .
 \end{equation*}
Then we have that $[\alpha]\in P^{p,q}(X)$ 
and $Q([\alpha],[\alpha])= 0.$
Applying Theorem A yields that $ [\alpha]=0.$ 
Hence,  $\tau$ is injective. \hfill $\square$

\smallskip

\noindent{\bf  Proof of Theorem C.}
Let $\phi:\  \H^{n-q,n-p}(X)\longrightarrow \H^{n-q+1,n-p+1}(X)$ be given by
\begin{equation*}
\phi([\alpha]):=[\omega_{n-p-q+1}]\wedge[\alpha],\qquad [\alpha]\in \H^{n-q,n-p}(X).
\end{equation*}
  Theorem B implies that $\dim P^{p,q}(X)=\dim \Ker \phi.$
On the other hand, by the classical Hodge theory (see \cite{bdi,gh,vo,we}) we know that $\phi$ is surjective. Hence,
\begin{eqnarray*}
\dim \Ker \phi & = & \dim \H^{n-q,n-p}(X)-\dim \H^{n-q+1,n-p+1}(X)\\
& = & \dim \H^{p,q}(X)-\dim \H^{p-1,q-1}(X).
\end{eqnarray*} 
Consequently,
\begin{equation}\label{dimen}
\dim P^{p,q}(X)+\dim \H^{p-1,q-1}(X)= \dim \H^{p,q}(X).
\end{equation}
On the other hand, it follows from Theorem B that the multiplication by 
$[\omega_{n-p-q+1}]$ is injective on $\H^{p-1,q-1}(X)$ and
\begin{equation}\label{inter}
P^{p,q}(X)\cap [\omega_{n-p-q+1}]\wedge\H^{p-1,q-1}(X)=\{0\}.
\end{equation}
Hence, the desired decomposition follows from (\ref{dimen}) and (\ref{inter}).  
\hfill $\square$

\section{Another version of the  Hodge--Riemann theorem}
In this section we describe the domain of validity 
of the mixed Hodge-Riemann theorem in the case where $p=q=1$. 
This problem is motivated by the dynamical study of holomorphic automorphisms
on compact K\"ahler manifolds.
An application of the mixed Hodge-Riemann theorem was given in the joint
work of the first author and Nessim Sibony \cite{ds1} (see also \cite{ds2}). 
In order to present  the  results  we need to introduce some notation.

\smallskip

Let $X$ be as usual a compact  K\"{a}hler manifold of dimension $n.$  
Define 
$$\H^{p,p}(X,\R):=\H^{p,p}(X)\cap \H^{2p}(X,\R).$$
Let $\K_p$ be the cone of all classes of smooth strictly positive $(p,p)$-forms 
in $\H^{p,p}(X,\R)$. This cone is open and satisfies $-\overline \K_p\cap
\overline \K_p =\{0\}$, where $\overline \K_p$ is the closure of $\K_p$. 
Each class in $\overline \K_p$ can be represented by a 
positive closed $(p,p)$-current.
The cone $\K_1$ is the {\it K\"ahler cone} of $X$.
Here, positivity of forms and currents of higher bidegree can be understood 
in the weak or strong sense.
We refer to \cite{de} for the basics on the theory of positive closed currents.

\smallskip

Fix a K\"ahler form $\omega$. Define $P^{p,q}(X)$ and $Q(\cdot,\cdot)$ 
as in (\ref{eq1.1}) and 
(\ref{eq1.2}) but 
for an arbitrary non-zero class $[\Omega]$ in $\overline\K_{n-p-q}$ and for 
$\omega_{n-p-q+1}:=\omega$.  The class $[\Omega]\wedge [\omega]$
does not vanish since it
can be represented by a non-zero positive closed current. 
Let $\K^{HR}_{n-p-q}$ be the cone of all classes $[\Omega]\in \K_{n-p-q}$ 
which satisfy the mixed Hodge-Riemann Theorem (Theorem A), 
that is,  $Q(\cdot,\cdot)$ is positive definite on $\P^{p,q}(X)$.

\smallskip

From now on we consider the case where $p=q=1$. 
The Poincar\'e duality implies that $\P^{1,1}(X)$ is a hyperplane of $\H^{1,1}(X)$ which depends 
continuously on $[\Omega]$. It follows by continuity that 
$\K^{HR}_{n-2}$ is an open cone in $\H^{n-2,n-2}(X,\R)$. Theorem A implies that one of  
the connected components 
of $\K^{HR}_{n-2}$ contains all the products of $(n-2)$ K\"ahler classes. Observe that 
$\P^{1,1}(X)$ does not intersect $\K_1$ since $[\Omega]$ is the class of a positive closed
current.

\smallskip

Let $\L_{n-2}$ be the set of all classes $[\Omega]$ in $\H^{n-2,n-2}(X,\R)$
such that the wedge product map $[\alpha]\mapsto
[\alpha]\wedge [\Omega]$ does not induce an isomorphism between 
$\H^{1,1}(X)$ and $\H^{n-1,n-1}(X)$.
Observe that 
$\L_{n-2}$ is an algebraic cone 
defined by a homogeneous polynomial of degree $\dim_\C\H^{1,1}(X)$.

\begin{prop}\label{otherversion} The cone
$\K^{HR}_{n-2}$ is a union of connected components of $\K_{n-2}\setminus \L_{n-2}$. 
In particular, it
does not depend on the K\"ahler form $\omega$. Moreover, if $[\Omega]$ is a class in 
$\overline\K^{HR}_{n-2}$ then $Q(\cdot,\cdot)$ is positive semi-definite on $\P^{1,1}(X)$ and 
for $c\in\P^{1,1}(X)$ we have
$Q(c,c)=0$  if and only
if $c\wedge[\Omega]=0$. 
\end{prop}
\begin{proof}
It is clear that $\L_{n-2}\cap \K^{HR}_{n-2}=\varnothing$.
Let $[\Omega]$ be a class in $\K_{n-2}$ which belongs to 
the boundary of $\K^{HR}_{n-2}$. We have to show that 
$[\Omega]\in \L_{n-2}$. By continuity, $Q(\cdot,\cdot)$ is positive semi-definite 
on $\P^{1,1}(X)$. Since 
$[\Omega]\not\in \K^{HR}_{n-2}$, there exists $c\in \P^{1,1}(X)$, $c\not=0$, such that $Q(c,c)=0$.
The Cauchy-Schwarz inequality implies that $Q(c,c')=0$ 
for every $c'$ in the hyperplane 
$\P^{1,1}(X)$. We have seen that $[\omega]$ does not belong to $\P^{1,1}(X)$.
On the other hand, 
$Q(c,[\omega])=0$ because $c\in\P^{1,1}(X)$. 
Consequently, $Q(c,\cdot)=0$. Therefore, the Poincar\'e duality
implies that $c\wedge [\Omega]=0$. Hence, $[\Omega]\in\L_{n-2}$. 
The first part of the proposition is proved.
We obtain the second part in the same way.
\end{proof}

\begin{rems}\rm
Observe that $Q(\cdot,\cdot)$ is positive definite on $\H^{2,0}(X)\oplus \H^{0,2}(X)$ for
$\Omega\in\K_{n-2}$. Then if $[\Omega]\in\K_{n-2}^{HR}$, the multiplication by $[\Omega]$
induces an isomorphism between $\H^2(X)$ and $\H^{n-2}(X)$.

\smallskip

Let $\widetilde \K_{n-2}$ be the cone of the classes in $\H^{n-2,n-2}(X,\R)$ 
of all positive closed 
currents of bidegree $(n-2,n-2)$.
This cone is convex and closed. Moreover, 
it contains $\K_{n-2}$ and satisfies $-\widetilde\K_{n-2}\cap 
\widetilde\K_{n-2}=\{0\}$.
Let $\widetilde\K^{HR}_{n-2}$ denote the cone of 
all classes in $\widetilde\K_{n-2}$ 
which satisfy the mixed Hodge-Riemann
theorem. Then Proposition \ref{otherversion} holds for 
$\widetilde\K^{HR}_{n-2}$. More precisely, $\widetilde\K^{HR}_{n-2}$ is
a union of connected components of $\widetilde\K_{n-2}\setminus \L_{n-2}$.
\end{rems}

The following type of results might be useful in the dynamical study of 
holomorphic automorphisms (see \cite{ds1,ds2})

\begin{cor}\label{limitcase} 
Let $c_1,\ldots,c_{n-2}$ be classes of  $\overline{\K}_1$ and let $c_{n-1}$ 
be a K\"ahler 
class. Then a class  $c$ in $\H^{1,1}(X)$ satisfies $c\wedge c_1\wedge \cdots 
\wedge c_{n-1}=0$ and
$c\wedge c\wedge c_1\wedge\cdots\wedge c_{n-2}  =0$ if and only if 
$c\wedge c_1\wedge\cdots\wedge
c_{n-2}  =0.$
\end{cor}
\begin{proof} 
Since each $c_i$ can be approximated by K\"ahler forms, 
Theorem A implies that $[\Omega]:=c_1\wedge\cdots\wedge c_{n-2}$
belongs to $\overline \K^{HR}_{n-2}$. Therefore, it is sufficient to apply
Proposition \ref{otherversion}.
\end{proof}

\end{document}